\documentclass{amsart}
\def \mce {{\mathcal E}}
\def \mcj {{\mathcal J}}
\def \length {\operatorname{length}}
\def \rank {\operatorname{rank}}

\renewcommand{\Re}{\operatorname{\rm Re}\nolimits}
\renewcommand{\Im}{\operatorname{\rm Im}\nolimits}

\def \Image {\operatorname{Image}}

\newtheorem{thm}{Theorem}[section]
\newtheorem{lemma}{Lemma}[section]
\newtheorem*{thmempty}{Theorem}
\newtheorem{prop}{Proposition}[section]
\newtheorem{cor}{Corollary}[section]
 
\newtheorem{defin}{Definition}[section]
\theoremstyle{definition}

\theoremstyle{remark}

\def\squarebox#1{\hbox to #1{\hfill\vbox to #1{\vfill}}}

\def \tr {\operatorname{tr}}
\def \dist {\operatorname{dist}}

\def \supp {\operatorname{supp}}

\def \comp {\operatorname{comp}}

\def \S1 {\mbox{S}^{1}}

\def \loc {\operatorname {loc}}

\def \Real {{\Bbb R}}

\def \Complex {\Bbb{C}}
\def \Natural {{\mathbb N}}
\def \Integers {{\mathbb Z}}

\def \rest {|}

   \title  [Resonances for manifolds with infinite cylindrical ends]
{Some Upper bounds on the number of resonances for manifolds
with infinite cylindrical ends}
   \author { T. Christiansen}
\thanks{Partially supported by an NSF grant.}
\thanks{e-mail address: {\tt tjc@math.missouri.edu}}

\begin{document}

\begin{abstract}

 We prove some sharp upper bounds on the number of resonances
associated with the Laplacian,
or Laplacian plus potential,
on a manifold with infinite cylidrical ends.
\end{abstract}

   \maketitle

The purpose of this note is to bound the 
number of resonances, poles of the meromorphic continuation
of the resolvent, associated to a manifold with infinite cylindrical
ends.  
These manifolds have an infinity which is in some sense one-dimensional,
even though the manifold is $n$-dimensional.  The bounds we obtain on the
resonances reflect this dichotomy- we obtain one bound like that for 
one-dimensional scattering, and other bounds of the type expected for
$n$-dimensional manifolds.  As part of our study of resonances, we relate 
poles of the resolvent to $L^2$ eigenvalues and to
poles of appropriately defined
``scattering matrices.''

\section{Introduction}
A smooth Riemannian 
manifold $X$ is said to be a manifold with cylindrical ends if
it 
can be decomposed as 
$$X= X_{\comp} \sqcup _{i=1}^{m_0}X_i,$$
where $X_i = [a_i, \infty)_{t} \times Y_i$, $(Y_i,g_i)$ is a compact 
Riemannian manifold, 
$X_{\comp}$ is a compact manifold with boundary $\sqcup_{1}^{m_0}Y_i$, and 
the metric on $X_i$ is $(dt)^2+g_i$.  We may also allow $X$ itself
to be a manifold with boundary, and then the ends may take the
form $[a_i, \infty)_{t} \times Y_i$, with $Y_i$ a 
smooth, compact manifold with boundary.
In this case, we require that the boundary of $X_c$ be compact and smooth
except for a finite number of corners corresponding to $a_i \times
\partial Y_i$.
An example of such a manifold is a waveguide,
a domain with 
smooth boundary in the plane, with one or
more infinite straight ends.  If we allow $X$ to have a boundary, we
will consider the Laplacian with Dirichlet or Neumann boundary conditions.

Let $\Delta_{Y_i}$ be the Laplacian on $(Y_i, g_i)$, and 
let $Y$ be the disjoint (and, if $m_0\geq 1$, disconnected) union
of the $Y_i$.  Let $\Delta _Y$ be the Laplacian on $Y$; that
is, if $f\in C^{\infty}(Y)$, then $(\Delta_Yf)_{\rest Y_i}
=\Delta_{Y_i}f_{\rest Y_i}$.  Let $\{ \sigma_j^{2}\}$, 
$\sigma_1^2\leq \sigma_2^2\leq \sigma_3^2\leq ...$ be 
the set of all eigenvalues of $\Delta_Y$, repeated according to
their multiplicity, and 
let $\nu_1^2< \nu^2_2 <\nu^2_3 <...$ be the {\em distinct}
eigenvalues of $\Delta_{Y}$.
Then the resolvent of the Laplacian $\Delta$ on $X$, or of $\Delta +V$, for
$V\in L_{\comp}^{\infty}(X)$ real-valued,
 has a meromorphic continuation to the 
Riemann surface $\hat{Z}$ on which $(z-\nu_j^2)^{1/2}$ is a 
single-valued function for all $j$  (\cite{guce, tapsit}).  Thus,
the resonances, poles of the meromorphic continuation of the 
resolvent, are associated to points in $\hat{Z}$.  The 
complicated nature of this Riemann surface makes more
difficult the question
of bounding the number of resonances, and makes necessary the restrictions
we place on $\{\nu_j^2\}$.

Part of our study of resonances is to understand the relationship
between poles of the resolvent and their multiplicities and poles
of a ``scattering matrix.''
For a manifold with cylindrical ends, there are several reasonable 
objects to call the scattering matrix.  One is an infinite dimensional
matrix.  In Section \ref{s:sm} we define this matrix and associated
finite-dimensional matrices which contain the information about poles
which we desire.  We show that if $z_0\in \hat{Z}$ is not a ramification
point of $\hat{Z}$ and is a pole of the resolvent, then there is an
associated $L^2$ eigenfunction with an appropriate expansion or there is a
pole of a ``scattering matrix'' at $z_0$.  
 We make this precise and address the
issue of multiplicities in Theorem \ref{thm:polesmult} and Proposition
\ref{p:notinspec}.  In \cite{a-p-v} such a relationship is noted, though
to our knowledge no proof appears in the literature, so we include it here
for completeness.

We will often require that there exists an $\alpha >0$ such that 
\begin{equation*}\tag{H1}
\nu_m^2-\nu^2_{m-1}\geq \alpha \nu_m  
\end{equation*}
for all sufficiently large $m$.  Examples of cross-sectional manifolds that 
satisfy such requirements are spheres, an interval with Dirichlet
or Neumann boundary conditions, and projective space.
We can easily construct a manifold that satisfies
(H1) with as many ends as we wish by 
taking all the cross sections $Y_i$ to be the same manifold from the 
previous list.  Other situations are possible as well, of course.

Let $P$ be the operator $\Delta$, the Laplacian, or $\Delta +V$,
for real-valued $V\in L^{\infty}_{\comp}(X)$.  Then, for 
$z\in \Complex \setminus [\nu^2_1, \infty)$, $R(z)=(P-z)^{-1}$ is bounded
on $L^2(X)$ except, perhaps, for a finite number of $z$.  It has a meromorphic
continuation to the Riemann surface $\hat{Z}$ described earlier.   
We bound the number of resonances in certain regions of $\hat{Z}$.
In doing so, we take the view that resonances near the physical sheet,
the sheet of $\hat{Z}$ on which the resolvent is bounded, are
more interesting, as they have greater physical relevance.  
Let 
$r_j(z)=(z-\nu_j^2)^{1/2}$.

Assuming the hypothesis (H1), we simplify the study of the resonances
somewhat and are then able to better bound them.

\begin{thm} \label{thm:ln}
Assume $X$ satisfies the hypothesis $(H1)$ and let $\beta <1$.
Then, in the connected components of $\{z'\in \hat{Z}: |r_m(z')|<\beta 
\sqrt{\alpha \nu_m}\}$ that meet the physical sheet, there are at most
${\mathcal O}(
m^{n-1})$ resonances.
\end{thm} 
We remark that 
we count all poles with their multiplicities
and that the poles of the resolvent include eigenvalues.
In \cite{ch-zw1, parnovski} 
there is an example of a family of manifolds that has
$\lim\inf_{\lambda\rightarrow \infty}N(\lambda)\lambda^{-n} >0$,
where $N(\lambda)$ is the number of eigenvalues of the Laplacian with norm
less than $\lambda^2$.  Since the cross-sectional manifolds in the
 example can be taken to be $n-1$-dimensional unit spheres,
this shows that the order appearing
in this theorem is optimal.
This type of bound is indicative of the $n$-dimensional
nature of the manifold.

In a simpler case, we can say a bit more.
\begin{thm}\label{thm:ppp}
Let $X= \Real \times Y$ and suppose $X$ satisfies (H1), and let 
$\rho>0$ be fixed.  Consider the operator $\Delta +V$,
for real-valued $V\in L^{\infty}_{\comp}(X)$.  
Then, on the connected components of $\{z\in \hat{Z}: |r_m(z)|<\rho\}$
that meet the physical sheet, the number of poles is bounded by
$C(1+m^{n-2})$.
\end{thm}
In case $V=V(s)$, $s\in \Real$, this theorem is easy to 
prove,
and this example shows that the bound is of optimal order.

Our proofs of these two theorems
 involve an adaptation of techniques developed by Melrose (\cite{pbdp}),
Zworski (\cite{zwrp,zwspb}),
 and Vodev (\cite{vodev})
 to bound the number of poles.   (See \cite{gu-zwncrs} and 
\cite{vodevsurvey} for further references and results.)
  It involves constructing
an approximation to the resolvent to find a holomorphic function whose
zeros include the poles of the resolvent.  Then we bound the function
and apply Jensen's theorem.  This
requires some knowledge of the function at a ``base point.''
For us, that will mean a lower bound.  Since we will be changing the
base point, we need some kind of uniform lower bound and that is 
different from these other applications of this technique.  In order
to do this, we will construct approximations of the resolvent 
especially well-suited to the regions where we work.

The following theorem does not require the hypothesis (H1), and 
its proof
uses a different technique.  
\begin{thm}\label{thm:fixedsheet}
Fix a sheet of $\hat{Z}$, and let $\{z_k\}$ be the resonances
of $P$ on this 
sheet.  Then 
$$\sum 
\frac{|\Im r_1(z_k)|}{|r_1(z_k)|^2}< \infty.$$
\end{thm}
This theorem is an analogue of what one finds in one-dimensional scattering
theory (see \cite{froese, zw1}), where the natural variable to consider is $\lambda = z^{1/2}$.

The problem of obtaining upper bounds on
the resonance-counting function has been widely
studied for Euclidean (e.g. \cite{froese,pbdp,vodev, zw1,zwspb}) and
hyperbolic scattering (e.g. \cite{ gu-zwncrs, gu-zwan,
muller, selberg}).  For a survey
and further references, see \cite{vodevsurvey} or \cite{zwsurvey}.

In this paper we use results of \cite{tapsit}, which studied the Laplacian
on compact manifolds with boundary and exact $b$-metrics.  Under a change
of variable, a special case of such manifolds is the class of manifolds
considered here (see also \cite{guce}).  The papers \cite{ch-zw1, parnovski}
independently obtained that the number of eigenvalues  less than $\lambda^2$
of the Laplacian grows at most like $\lambda^n$.  The existence of
eigenvalues or complex resonances has been studied in, for 
example, \cite{a-d, a-p-v, b-g-r-s, ch-zw1, 
d-p, d-e-m, parnovski} and references.

In finishing the paper, the author received a copy of 
\cite{edward}.  There Edward obtains a result similar to our Theorem
\ref{thm:ln} for the Laplacian on waveguides, that is, domains in
the plane which outside of a compact set coincide with $(-\infty, \infty)
\times \pi$, and thus fall in the category of manifolds which we consider.
The waveguides satisfy hypothesis (H1) and their $\nu_m^2$ are 
very regularly distributed for either the Dirichlet or Neumann Laplacian.
Use the metric induced on $\hat{Z}$ by the pull-back of the metric 
on $\Complex$ to define distance on 
$\hat{Z}$.  Then summing over the $\nu_m^2 \leq r$ we can obtain from
our Theorem \ref{thm:ln} an upper bound 
\begin{multline*}
\# \{z_j: z_j \; \text{is a resonance of }\; \Delta, \; \\
\dist(z_j,\text{ physical sheet})<c_1\sqrt{|\pi(z_j)|},\;  |\pi(z_j)|<r \}
  =
 \mathcal{O}(r)
\end{multline*} for some $c_1>0$,
where $\pi:\hat{Z}\rightarrow Z$ is projection.  
This may be compared to Theorem 2 of 
\cite{edward}, where for $c_1=1/2$, Edward obtains a bound 
$\mathcal{O}(r^{3+\epsilon})$, any $\epsilon >0$.

{\bf Remark.} We note that the ``black-box'' formalism of \cite{sj-zw}
can be adapted to this situation.  Thus we could replace $P$ by 
$\tilde{P}$, a more general
self-adjoint, compactly supported perturbation of the 
Laplacian than the ones considered here.  If $\tilde{P}$ satisfies 
the assumptions of \cite{sj-zw}, properly interpreted for this setting,
and is bounded below,
then Theorem \ref{thm:fixedsheet} will hold for $\tilde{P}$.  Let
$\tilde{P}^{\#}=\tilde{P}_{\rest \{t<\max(a_i)\}},$ with Dirichlet 
boundary conditions at $\{t=\max(a_i)\}.$  Let 
$$N_{\tilde{P}^{\#}}(\lambda)=\# \{ \lambda^2_j: \lambda_j^2 \; \text{is 
an eigenvalue of }\; \tilde{P}^{\#},\; \lambda_j^2 \leq \lambda^2\}.$$
If, in addition to the other assumptions,
$$N_{\tilde{P}^{\#}}(\lambda+a)-N_{\tilde{P}^{\#}}(\lambda-a)=
\mathcal{O}_a(\lambda^q),$$ $q \geq n-1$, for any constant $a$, then an
analogue of Theorem \ref{thm:bbn} holds, with $\mathcal{O}_{\beta}(m^{n-1})$
replaced by $\mathcal{O}_{\beta}(m^q)$.

{\bf Acknowledgments.} It is a pleasure to thank Maciej Zworski for 
suggesting the problem of counting resonances in this setting.  I am 
grateful to him and to Dan Edidin for helpful discussions.  Thanks to
Julian Edward for helpful comments.

\section{Preliminaries}

Let $r_j(z)=(z-\nu_j^2)^{1/2}$ and 
identify the physical sheet of 
$\hat{Z}$ as being the part of $\hat{Z}$ on which $\Im r_j(z)>0$ for
all $j$ and all 
$z$ and on which $R(z)$ is bounded.  Other sheets will be identified,
when necessary, by indicating for which values of $j$ $\Im r_j(z)<0$. 
Each sheet can be identified with  $\Complex \setminus [\nu_1^2,\infty)$.
With this language, there are points in $\hat{Z}$ which belong to no sheet
but which belong to the boundary of the closure of two sheets, and the 
ramification points, which correspond to $\{\nu_j^2\}$ and belong to 
the closure of four sheets (except for ramification points
corresponding to $\nu_1^2$).  We note that sheets that meet the physical sheet
are characterized by the existence of a $J\in \Natural $ such that
$$\Im r_j(z)<0 \; \text{for all } z \;\text{on that sheet if and only if} \;
j \leq J.$$

Let $\{\phi_j\}$ be an orthonormal set of eigenfunctions of $\Delta_Y$
associated with $\{\sigma_j^2\}$.

We shall often abuse notation and use $y$ for a coordinate in either
$Y$ or $Y_i$.  On an end, we use the coordinates $(t, y)$, with $t\in (a_i,
\infty)$ and $y\in Y_i$.

We define an operator on $[0,\infty)_t \times Y_y$.
Using the same notation for an operator and its Schwartz kernel, let
$$R_{ej}(z)=\frac{i}{2(z-\sigma_j^2)^{1/2}}(e^{i|t-t'|(z-\sigma_j^2)^{1/2}}
-e^{i|t+t'|(z-\sigma_j^2)^{1/2}})\phi_j(y)\overline{\phi_j}(y'),$$
where we take $\Im (z-\sigma_j^2)^{1/2}>0$ for $z\in \Complex
\setminus [\sigma_1^2,\infty)$ corresponding to the physical sheet of
$\hat{Z}$.
Then, for $z\in \Complex \setminus [\sigma_1^2,\infty)$,
\begin{equation}\label{eq:RY}
R_Y(z)=\sum_j R_{ej}(z)
\end{equation}
is the resolvent of the Dirichlet Laplacian on $[0,\infty)\times Y$.
As an operator from $L^2_{\comp}([0,\infty)\times Y)$ to 
$H^2_{\loc}
([0,\infty)\times Y)$ it has a holomorphic continuation to the 
Riemann surface $\hat{Z}$ and we will see that it is this that determines
the surface to which the resolvent of $P$ has a meromorphic
continuation.

In general, we shall use $z$ to stand for a point in $\hat{Z}$ and 
$\pi(z)$ to represent its projection to $\Complex$.  For 
$w\in \Real^m$, $\langle w \rangle =(1+|w|^2)^{1/2}$.

\section{The Scattering Matrix and Poles of the Resolvent}\label{s:sm}

There are several reasonable definitions of the scattering matrix
in this setting.  We recall some results of \cite{ace, tapsit}.

Recall that 
$r_k(z)=(z-\nu_k^2)^{1/2}$,
and let $\tilde{r}_j(z)=(z-\sigma_j^2)^{1/2}$, with $\Im \tilde{r_j}>0$ when
$z$ is in the physical sheet and $\tilde{r}_j(z)=r_k(z)$ if 
$\sigma_j^2=\nu_k^2$.  In addition, fix a coordinate $t$ on the ends, so that $t=0$ lies on the ends.  Then, for all but a finite number of $z$ 
in the physical space, there are functions $\Phi_j(z,p)$ with
\begin{equation}\label{eq:phij1}
(P-\pi(z))\Phi_j(z,p)=0
\end{equation}
and, on the ends, 
\begin{equation}\label{eq:phij2}
\Phi_j(z,t,y)=e^{-i\tilde{r}_j(z)t}\phi_j(y)+\sum S_{mj}(z)
e^{i\tilde{r}_m(z)t}\phi_m(y).
\end{equation}
 The $\Phi_j$ have a meromorphic continuation
to all of $\hat{Z}$ and thus, so do the $S_{mj}$.  The $S_{mj}(z)$
depend on the choice of the coordinate $t$ in a fairly straight-forward
way (see \cite{ch-zw1}).  This dependence is not important here
as it does not change the location of the scattering poles, so we 
ignore it but we do consider the coordinate $t$ to be fixed throughout,
and chosen so that $\{p\in X: t=0\}\subset \sqcup X_i$. 

There are several reasonable choices of objects to 
call the scattering matrix.  One possibility is the 
infinite matrix of the $S_{ij}(z)$, as in \cite{tapsit}.  Another,
which is well-defined for $z$ on the boundary of the physical sheet,
is a normalized, finite-dimensional matrix of the $S_{ij}(z)$, where the 
dimension changes as $z$ crosses a $\nu_j^2$.  This is used in \cite{ace}
and has the advantage of being unitary (though we note that
the variable used in \cite{ace} is $\lambda = z^{1/2}$).  Here, 
however, this is unnecessarily complicated as it requires the introduction
of $(z-\nu_i^2)^{1/4}$.  

We shall work with finite-dimensional matrices, of the form 
$$(S_{ij}(z))_{i,j\in \mce}$$
for some set $\mce \subset \Natural$, 
where $\mce$ is chosen, depending on $z$, to be most helpful
for our purposes.

\subsection{Definitions and Some Properties}
In order to describe the poles and their multiplicities
 of the resolvent and matrices of the $S_{ij}(z)$, we introduce some 
notation.

\begin{defin} We define the multiplicity of a pole of the resolvent
$R(z)$ at $z_0$ to be
$$m_{z_0}(R)=\dim \Image \Xi_{z_0}(R),$$
where $\Xi_{z_0}(R)$ is the singular part of
$R$ at the point $z_0.$
\end{defin}
It follows as in \cite[Lemma 2.4]{gu-zwan} that, if $z_0$ is not a ramification
point of $\hat{Z}$, then 
$m_{z_0}(R)$ is also the rank of the residue of $R$ at $z_0$.  

In order to define the multiplicity of the pole 
of a matrix, we shall use the following lemma.
  Though it may be well known, we 
include it and a proof for the convenience of the reader.
\begin{lemma}\label{l:matform}
Suppose $A(z)$ is a $d\times d$-dimensional meromorphic matrix, 
invertible for some value of $z$.  
Then, near $z_0$, it can be put into
the form 
\begin{equation*}
A(z) = E(z)\left(\sum_{j=1}^p(z-z_0)^{-k_j}P_j + \sum _{j=p+1}^{p'}
(z-z_0)^{l_j}P_j+P_0\right)F(z)
\end{equation*}
where $E(z)$, $F(z)$, and their inverses are holomorphic near $z_0$, and
$P_iP_j=\delta_{ij}P_i$, $\tr P_0= d-p'$, $\tr P_i=1$, $i=1,...,p'$.
The $k_j$ and $l_j$ are, up to rearrangement, uniquely determined.
\end{lemma}  
\begin{proof}
We outline a proof.  Without loss of generality we may assume that 
$z_0=0$.  

First we show the existence of such a decomposition.  Choose $k$ such that 
$z^kA(z)=B(z)$ is holomorphic at $0$.  Now the proof of the existence
of such a decomposition follows much like a proof from \cite[Section VI.2.4]
{gantmacher}.  Let $B(z)=(b_{ij}(z))$.  Choose an element $b_{ij}(z)$ that
vanishes to the lowest order at $0$, and by permuting rows and columns
make this element $b_{11}(z)$.  By subtracting from the $k$th row the first 
row multiplied by $b_{k1}(z)/b_{11}(z)$
(which is holomorphic near $z=0$), we can make all the entries in the first 
row, other than the first one, zero.  Similar column operations reduce $B(z)$
to the form
\begin{equation*}
\left( \begin{array}{cccc}
b_{11}(z) & 0 & \cdot \cdot \cdot & 0\\
0 & c_{22}(z) & \cdot \cdot \cdot & c_{2d}(z)\\
\vdots & \vdots & \cdots & \vdots \\
0 & c_{d2}(z) & \cdot \cdot \cdot & c_{dd}(z) 
\end{array}
\right)_.
\end{equation*}
Repeating the procedure on the $(d-1)\times(d-1) $
dimensional matrix of the $c_{ij}$, we obtain a matrix of the form
\begin{equation*}
\left( \begin{array}{ccccc}
b_{11}(z) & 0 & 0 &\cdot \cdot \cdot & 0\\
0 & c_{22}(z) & 0 & \cdots & 0\\
0 & 0& \alpha_{33}(z) & \cdots & \alpha_{3d}(z)\\
\vdots & \vdots & \vdots & \hdots &  \vdots \\
0 &0 & \alpha_{d3}(z) & \cdot \cdot \cdot & \alpha_{dd}(z) 
\end{array}
\right)_.
\end{equation*}
Repeating this procedure a finite number of times,
 we obtain a matrix with only the diagonal
entries, $h_{ii}(z)$, nonzero.  Since $h_{ii}(z)=z^{m_i}g_i(z)$, with
$g_i(0)\not = 0$ and 
$g_i(z)$ holomorphic near $0$, by multiplying by a diagonal matrix
with nonzero entries $1/g_i(z)$ we obtain a diagonal
matrix with diagonal entries of the form $z^{m_i}$.  Finally, multiplying
the whole thing by $z^{-k}I$ we obtain a matrix equivalent to $A(z)$ near
$z=0$.  We note that the construction ensures that $E$ and $F$ are holomorphic
near $0$.  Moreover, since each of $E$ and $F$ is the product of elementary
matrices and one diagonal matrix (holomorphic near
$0$, with nonzero determinant), $E$ and $F$ are both invertible near $0$.

To see the uniqueness of the $k_i$, $l_i$, we prove it for $B(z)=z^{k}A(z)$,
where $B$ is holomorphic at $0$, and
and use the straight-forward relationship between $A$ and $B$.
Suppose there are two such decompositions: $B(z)=E_1(z)D_1(z)F_1(z)
=E_2(z)D_2(z)F_2(z)$, where $E_i(z)$ and $F_i(z)$ are as in the statements
of the lemma and the $D_i$ are diagonal matrices with nonzero entries
$d_{mm,i}=z^{l_{m,i}}$.  By row and column operations we can make  
$0\leq l_{1,i}\leq l_{2,i}\leq ...\leq l_{d,i}$.  We have
\begin{equation}\label{eq:me1}
E(z)D_1(z)=D_2(z)F(z)
\end{equation}
for new matrices $E(z),F(z),$ holomorphic and invertible 
 near $0$.  Thus it is
easy to see that 
$$\rank(D_1(0))=\rank(D_2(0))\equiv r_0.$$
Using (\ref{eq:me1}) and the fact that 
$l_{j,i}=0$ if and only if $j\leq r_0$, we obtain that, if
$E(z)=(e_{ij}(z))$,  $F(z)=(f_{ij}(z))$,
\begin{equation}\label{eq:ef0}
e_{ij}(0)=0 = f_{ji}(0)\;\text{ if}\; j\leq r_0\; \text{ and} \;i>r_0.
\end{equation}

We finish the proof of the uniqueness by induction.  Let $r_{j,i}=\rank(D^{(j)}
_i(0)), $ and notice it suffices to prove that 
$r_{j,1}=r_{j,2}$ for all $j$, and that we have
$\sum_{j}r_{j,i}=d$.  Suppose we have shown that 
$r_{q,1}=r_{q,2}\equiv r_q$ for $q\leq N$, and that, if 
$R_N=r_0+r_1+\cdots + r_N$, 
\begin{equation}\label{eq:efN}
e_{ij}(0)=0=f_{ji}(0)\; \text{for}\; j \leq R_N, \; i>R_N.
\end{equation}
If $v=(v_1,v_2,...,v_d)^t$, let $\Pi_{s}v=(0,...,0, v_{s+1},...,v_d)^t$
for $s\in \Natural$, $s\leq d$.  Then 
$$\rank(\Pi_{R_N}E(0)\Pi_{R_N})=d-R_N=\rank(\Pi_{R_N}F(0)\Pi_{R_N}),$$
using (\ref{eq:efN})
and the fact that $E(0)$ and $F(0)$ are invertible.  Since,
for $j\leq N$, $D^{(j)}_1(0)\Pi_{R_N}=0$ and $\Pi_{R_N}D^{(j)}_2(0)=0$,
we get
\begin{multline*}
\rank\left(\frac{d^{N+1}}{dz^{N+1}}(\Pi_{R_N}E(z)D_1(z)\Pi_{R_N})_{\rest z=0}
\right)
= \rank \left(\Pi_{R_N}E(0)D_1^{(N+1)}(0)\Pi_{R_N}\right) \\
= \rank\left(D_1^{(N+1)}(0)\right) 
= 
\rank\left(\frac{d^{N+1}}{dz^{N+1}}(\Pi_{R_N}D_2(z)F(z)\Pi_{R_N})_{\rest z=0}
\right) \\
=\rank \left(D_2^{(N+1)}(0)\right)
\end{multline*}
and thus $r_{1,N+1}=r_{2,N+1}\equiv r_{N+1}.$  If $R_{N}+r_{N+1}<d$, from 
(\ref{eq:me1}) we obtain
$$e_{ij}(0)=0=f_{ji}(0)\; \text{for}\; j\leq R_N+r_{N+1}, \;
i>R_N+r_{N+1}.$$
The induction need only continue until $R_{N+1}=d$, finishing the proof.
\end{proof}

\begin{defin}  Let 
$A(z)$ be a meromorphic matrix, invertible for some 
values of $z$, and let $k_1, k_2,...k_p$, $l_{p+1},l_{p+2},...l_{p'}$ be as in 
Lemma \ref{l:matform}.  Set
$$\mu m_{z_0}(A)=\sum_{j=1}^pk_j,$$
the ``maximal multiplicity'' of the pole of $A$ at $z_0$.  Set
$$\mu d_{z_0}(A) = \sum_{j=1}^pk_j - \sum _{j=p+1}^{p'}l_j, $$
the ``determinantal multiplicity'' of the pole of $A$ at $z_0$.  
\end{defin}We 
note that 
$$\mu d_{z_0}(A)= \min \{j \in \Integers: (z-z_0)^j \det A(z) \text{ 
is regular at } z_0\}.$$
Each of these measures of multiplicity will be useful.

Suppose $\mce\subset \Natural$ is a finite subset.  Let
$$\tilde{\mce}=\{ j\in \Natural : \sigma_j^2=\nu_l^2\;
 \text{for some } \; l\in \mce\}.
$$
Define the matrix 
$$S_{\mce}(z)=(S_{ij}(z))_{i,j\in \tilde{\mce}}.$$
For $z\in \hat{Z}$, let 
$$
\mce_z=\{ j\in \Natural : \Im r_j(z) \leq 0\}$$
and $$\mcj_z= \{ j\in \Natural : j \leq \max \mce_z \}.$$

If $\mce \subset \Natural$ is a finite set, define 
$w_{\mce}:\hat{Z}\rightarrow \hat{Z}$ as follows.
To $z$ we may associate the set of square roots $\{r_j(z)\}$.  Then 
$w_{\mce}(z)$ may
be determined by saying it it the element of $\hat{Z}$  
associated to the set $\{ r_j(w_{\mce}(z))\}$, with
$$r_j(w_{\mce}(z))=\left\{ \begin{array}{ll} -r_j(z),\; &
 \text {if}\; j \in \mce \\
r_j(z),& \text {if}\; j\not \in \mce.
\end{array} \right.
$$
For $J\in \Natural$, define a similar operator 
$w_J:\hat{Z}\rightarrow \hat{Z}$, where 
$$r_j(w_J(z))=\left\{ \begin{array}{ll} -r_j(z),\; & \text {if}\; j \leq J\\
r_j(z),& \text {if}\; j>J.
\end{array}
\right.
$$
This operator appears in the relations (\ref{eq:drphi}) and (\ref{eq:ri}).

The following Proposition generalizes similar results of \cite{tapsit} and
\cite{ace}.
\begin{prop}\label{prop:sinverse}
Let $\mce \subset \Natural$ be a finite set.
For $j\in \tilde{\mce}$, 
$$\Phi_j(w_{\mce}(z))=\sum _{k \in \tilde{\mce}}S_{kj}(w_{\mce}(z))\Phi_k(z)$$
and $$(S_{\mce}(z))^{-1}=S_{\mce}(w_{\mce}(z)).$$
\end{prop}
\begin{proof}
For $j\in \tilde{\mce}$, consider
$$\Psi_j(z) =\Phi_j(w_{\mce}(z))-\sum_{k\in \tilde{\mce}}S_{kj}(w_{\mce}(z))
\Phi_k(z).$$
Then
$$(P-\pi(z))\Psi_j(z)=0$$
and, on the ends, $\Psi_j$ has an expansion
\begin{multline}
\Psi_j(z,t,y)=
e^{i\tilde{r}_j(z)t}\phi_j (y) 
-\sum_{k,l\in \tilde{\mce}}
S_{kj}(w_{\mce}(z))S_{lk}(z)e^{i\tilde{r}_l(z)t}\phi_l(y)\\
+\sum_{l\not \in \tilde{\mce}}\left(S_{lj}(w_{\mce}(z))- 
\sum_{k\in \tilde{\mce}}S_{kj}(w_{\mce}(z))S_{lk}(z)\right)
 e^{i\tilde{r}_l(z)t}\phi_l(y).
\end{multline}
Since, for $z$ on the physical sheet, $\Im \tilde{r}_l(z)>0$ for 
all $l$, $\Psi_j(z)\in L^2(X)$ there and thus $\Psi(z)\equiv 0$ for
all $z$ in the physical space.  By analytic continuation, $\Psi_j(z)\equiv 0$
for all
$z \in \hat{Z}$.  This proves the first part of the Proposition.  It also
shows that, for $j,l\in \tilde{\mce}$,
$$\sum_{k \in \tilde{\mce}}S_{kj}(w_{\mce}(z))S_{lk}(z)=\delta_{jl};
$$ that is, $(S_{\mce}(z))^{-1}=S_{\mce}(w_{\mce}(z)).$
\end{proof}

\subsection{Relation between Poles of the Resolvent and Poles of 
$S_{\mcj{z_0}}(z)$}

First, we recall some results of \cite[Sections 6.7, 6.8]{tapsit} on the 
nature of poles of the resolvent on the boundary of the physical space.
\begin{prop}
Suppose $z$ lies on the boundary of the physical space.  Then, if $z$ 
is not a ramification point of $\hat{Z}$, the multiplicity of $z$ 
as a pole of the resolvent is equal to the dimension of the $L^2$ 
null space of $P-\pi(z)$.  If $z$ is a ramification 
point, then the resolvent has a double pole at $z$, with coefficient
the projection onto the $L^2$ null space of $P-\pi(z)$, if there is
any.  The residue has rank equal to the dimension of 
$$\{f: (P-\pi(z))f=0, \; f\not \in L^2(X),\; \frac{\partial}{\partial t}f 
\in L^2(X) \}.$$
\end{prop}

 Since
\begin{equation}\label{eq:cge}
\Phi_j(z,p)=\chi(t)e^{-i\tilde{r}_j(z)t}\phi_j
-(P-z)^{-1}(P-\pi(z))
\chi(t)e^{-i\tilde{r}_j(z)t}\phi_j,
\end{equation}
where $\chi(t)\in C^{\infty}(\Real)$ is supported in $t>\max(a_i)$
and is one in a neighborhood of infinity,
it is clear by their definition that the $S_{ij}(z)$ cannot have a pole unless
$R(z)$ has a pole.  
More can be said, and we begin our study of the relationship between the 
poles of the resolvent and 
poles of the scattering matrix with the following
\begin{thm}\label{thm:polesmult}
Suppose $z_0\in \hat{Z}$, and $z_0$ is not a ramification point.  Then
\begin{multline*}
m_{z_0}(R)
\\=\mu m_{z_0}(S_{\mcj_{z_0}}) + 
\dim \{f\in L^2(X): (P-\pi(z))f=0, f\sim \sum_{j\not\in  \tilde{\mcj}_{z_0}} c_j e^{ir_j(z)t}\phi_j\}
\end{multline*}
where the expansion must be valid on any end.
\end{thm}
In proving this theorem, we will use some techniques from 
\cite[Section 2]{gu-zwan}.

We will call 
\begin{equation}\label{eq:typeI}
\{f: (P-\pi(z))f=0, f\sim \sum_{j\not\in  \tilde{\mcj}_{z_0}} c_j e^{ir_j(z)t}\phi_j\}
\end{equation}
the set of 
eigenfunctions of type I.  This depends on $z_0$ of course, and we 
will note the dependence in cases of possible confusion.

We first show that
\begin{lemma}\label{l:mumsmaller}
 If $z_0$ is not a ramification point of $\hat{Z}$, then
\begin{multline*}
\mu m_{z_0}(S_{\mcj_{z_0}}) \\
\leq m_{z_0}(R)-\dim \{f\in L^2(X): (P-\pi(z))f=0, f\sim \sum_{j\not\in  \tilde{\mcj}_{z_0}} c_j e^{ir_j(z)t}\phi_j\}.
\end{multline*}
\end{lemma}
\begin{proof}

Just as in Lemma 2.4 of \cite{gu-zwan}, whose notation we use, one can prove that if 
$R(z)$ has a pole at $z_0$, for $z_0$ not a ramification point of $\hat{Z}$,
then near $z_0$
\begin{equation}\label{eq:Rpoles}
R(z)=\sum_{k=1}^p \frac{A_k(z_0)}{(z-z_0)^k}+H(z_0,z)
\end{equation}
where $$A_{k}(z_0)=\sum_{l,m=1}^q a_k^{lm}(z_0) \varphi_l \otimes
\varphi_m$$
and $$(\varphi_l \otimes \varphi_m) f (p)=\varphi_l(p) \int_{p'\in X}
\varphi_m(p')f(p').$$
As in \cite{gu-zwan}, 
we have $$(P-\pi(z_0))A_k(z_0)=A_{k+1}(z_0)=A_k(z_0)(P-\pi(z_0)).$$
If $a_k(z_0)=(a_k^{lm}(z_0))_{1 \leq l,m \leq q}$, then 
$a_1(z_0)$ is symmetric with rank $q$, $d(z_0)=a_1(z_0)^{-1}a_2(z_0)$
is nilpotent, and $a_k(z_0)=a_1(z_0)d^{k-1}(z_0),$ $k>1$.  Moreover, 
$\varphi_l$ has an expansion on the ends of the form 
\begin{equation}\label{eq:phiexp}
\varphi_l(t,y)_{\rest t>\max a_i} 
 =\sum_j \sum _{m\leq m_l} c_{ljm} t^m e^{i\tilde{r}_j(z_0)t}
\phi_j(y)
\end{equation}
(\cite{tapsit})
and the $\varphi_l$ are linearly independent.  Let $B(z)$ be the 
matrix $B(z)=(b^{lm}(z))_{l,m \leq p}$ with 
$$b^{lm}(z) = \sum _{k=1}^p \frac{a_k^{lm}(z_0)}{(z-z_0)^{k}}.$$
Then $B(z)$ can be written, as in Proposition 2.11 of
\cite{gu-zwan}, as
\begin{equation}\label{eq:nicebrep}
B(z)= E^{\#}(z) \left[ \sum_{j=1}^{p'}(z-z_0)^{-k_j}P_j+P_0\right] F^{\#}(z)
\end{equation}
where $E^{\#}(z)$, $F^{\#}(z)$ and their inverses are holomorphic
 near $z_0$, 
$P_iP_j=\delta_{ij}P_i$, $\tr P_i=1$, $i\not = 0$, $\tr P_0=q-p'$,
and $k_1+k_2+...+k_{p'}=q$.

First assume that there are no type I eigenfunctions of $P$ with eigenvalue
$\pi(z_0)$.
We recall that we can construct the 
generalized eigenfunctions $\Phi_j$,
$j\in \tilde{\mcj}_{z_0}$, as indicated in (\ref{eq:cge}).  We 
then obtain the entries of $S_{\mcj_{z_0}}(z)$ by first restricting 
$\Phi_j$ to $t=0$ and then extracting the 
coefficients of $\phi_k$.  

This, taken with the representation 
(\ref{eq:nicebrep}), shows that the singular part of $S_{\mcj_{z_0}}(z)$ near
$z_0$ will be given by 
$$E^{\flat}(z)\left[ \sum(z-z_0)^{-k_j}P_j+P_0\right]F^{\flat}(z)
$$ 
where $E^{\flat}(z)$, $F^{\flat}(z)$ are holomorphic matrices 
near $z_0$ which may 
not be invertible, or even square, and this proves the lemma in this 
special case.

To handle the case where $\pi(z_0)$ is an eigenvalue of $P$ with
\begin{equation}
n_e(z_0)=\dim \{f\in L^2(X): (P-\pi(z_0))f=0, f\sim \sum_{j\not\in  
\tilde{\mcj}_{z_0}} c_j e^{ir_j(z)t}\phi_j\} >0,
\end{equation}
we need to be a bit more careful.  We first remark that for a general
$z\in \hat{Z}$, $R(z)\psi$ is well-defined if
$$\psi \in e^{-\langle t\rangle\max (0, -\Im r_j(z))}L^2(X).$$
This can be seen from the construction of the analytic continuation
of the resolvent in \cite{tapsit}.
Now suppose $\psi$ is an $L^2$ eigenfunction of $P$.  Then 
it is exponentially decreasing, and, for $J\in \Natural$, if 
$z$ and $w_J(z)$ both lie on the boundary of the physical space,
\begin{equation}\label{eq:drphi}\left( R(z)-R(w_J(z))\right)\psi =0.
\end{equation}
By analytic continuation, if $z$, $w_J(z)$ are such that 
$$\psi \in 
e^{-\langle t \rangle \max (0, -\Im r_j(z), -\Im r_j(w_J(z)))}L^2(X),$$
then (\ref{eq:drphi}) holds.  By repeatedly applying (\ref{eq:drphi})
and using our knowledge of the structure of the 
resolvent on the closure of the physical space, we obtain that if
$$\psi\in \{f\in L^2(X): (P-\pi(z_0))f=0, f\sim 
\sum_{j\not\in  \tilde{\mcj}_{z_0}} c_j e^{ir_j(z)t}\phi_j\},$$
then there is a neighborhood of $z_0$ so that in this neighborhood
$$R(z)= \frac{1}{(\pi(z_0)-\pi(z))\|\psi\|^2}\psi \otimes \overline{\psi}
 + B_1(z,z_0,\psi)$$
with $B_1(z,z_0,\psi)\psi =0.$  That
is, each $L^2$ eigenfunction with an expansion at infinity of this
type contributes to the singularities of $R(z)$ in the expected way.
  Using the fact that if $j\in \mcj_{z_0}$
$$\int _{X} \overline{\psi}
 (P-\pi(z_0))\chi(t)e^{-i\tilde{r}_j(z_0) t}\phi_j(y) = 0,$$
we see that the singularities of $R(z)$ at $z_0$ corresponding to $L^2$
eigenfunctions of type I
do not contribute to the singularities of $S_{\mcj_{z_0}}(z)$ at $z_0$.  
After making this observation, the proof for the case $n_e(z_0)>0$
follows just as in the case $n_e(z_0)=0$.
\end{proof}

The proof of Theorem \ref{thm:polesmult} will be completed by 
\begin{lemma}\label{l:mugreater}
If $z_0$ is not a ramification point of $\hat{Z}$, then
\begin{multline*} \mu m_{z_0}(S_{\mcj_{z_0}})+
\dim \{f\in L^2(X): (P-\pi(z))f=0, f\sim 
\sum_{j\not\in  \tilde{\mcj}_{z_0}} c_j e^{ir_j(z)t}\phi_j\} \\
\geq m_{z_0}(R).
\end{multline*}
\end{lemma}
\begin{proof}

We use, for $J\in \Natural$,
\begin{equation}\label{eq:ri}
R(z)-R(w_J(z))
= \frac{i}{2} \sum_{\sigma_j^2 \leq \nu_J^2}
\Phi_j(z)\otimes \Phi_j(w_J(z))\frac{1}{\tilde{r}_j(z)}.
\end{equation}
Equation (\ref{eq:ri}) holds, by Stone's formula, \cite[Section 6.8]{tapsit},
and \cite[Section 2.2]{ace},
 for $z$ on the boundary of the physical
space, with $\nu_{J}^2 < \pi(z) < \nu_{J+1}^2$, and then holds on the
rest of $\hat{Z}$ by analytic continuation.
Using Proposition \ref{prop:sinverse},
we may write (\ref{eq:ri}) as
\begin{equation}
R(z)-R(w_J(z))
= \frac{i}{2} \sum_{\sigma_j^2 \leq \nu_J^2}\sum_{\sigma_m^2 \leq \nu_J^2}
 S_{mj}(z)\Phi_m(w_J(z))\otimes \Phi_j(w_J(z))\frac{1}{\tilde{r}_j(z)}.
\end{equation}  

Recalling equations (\ref{eq:Rpoles}) and (\ref{eq:phiexp}), we see that
if $\psi$ is in the image of the singular part of $R(z)$ at $z_0$, then
$\psi$ is a linear combination of
eigenfunctions of type I (see (\ref{eq:typeI}))
and 
\begin{multline}\label{eq:typeII}
\{ f : (P -\pi(z_0))^k=0\; \text{for some }k \in \Natural; \; 
f \sim  \sum_j \sum _{m\leq m_l} b_{jm} t^m e^{i\tilde{r}_j(z_0)t}
\phi_j(y), \\ b_{jm} \not = 0\; \text{for some}\; j \in \tilde{J}_{z_0},
\text{some}\; m
\} .
\end{multline} 
We call functions of the form (\ref{eq:typeII}) type II.  Now
take $J=\max{\mce}_{z_0}$.
If $g $ is in the image of the singular
parts of $R(z)$ at both
$z_0$ and at $w_J(z_0)$, then it must be of type I.  

It is the appearance of the type II functions which we must understand.
Suppose $g$ is in the image of the residue of $R(z)$ at $z_0$ and is of
type II.  Then, since it is not in the image of the residue of 
$R(z)$ at $w_J(z_0)$, it must be in the image of the residue of
\begin{equation}\label{eq:diff}
\sum_{\sigma_j^2 \leq \nu_J^2}\sum_{\sigma_m^2 \leq \nu_J^2}
 S_{mj}(z)\Phi_m(w_J(z))\otimes \Phi_j(w_J(z))\frac{1}{\tilde{r}_j(z)}
\end{equation}
at $z=z_0$.  First assume that $\Phi_j(w_J(z))$, $j\in \tilde{J}_{z_0}$
has no poles at $z_0$.  Recalling that we may write $S_{\mcj_{z_0}}(z)$
near $z_0$ as 
in Lemma \ref{l:matform}, we can write (\ref{eq:diff}) as
\begin{equation}\label{eq:rewritediff}
\sum_{j=1}^p (z-z_0)^{-k_j}\Psi_j^1(z)\otimes \Psi_j^2(z)
+ H(z,z_0)
\end{equation}
where $k_1+...+k_p= \mu m_{z_0}(S_{\mcj_{z_0}}(z))$, and $\Psi^1_j(z)$, 
$\Psi_j^2(z)$,
and $H(z,z_0)$
are holomorphic near $z_0$.  Then it is easy to see that the rank of the 
image of the residue of (\ref{eq:rewritediff}) cannot exceed 
$\mu m _{z_0}(S_{\mcj_{z_0}}(z))$, and the total rank of the residue of
$R(z)$ at $z_0$ cannot exceed
$$\mu m_{z_0}(S_{{\mcj}_{z_0}})+
\dim \{f\in L^2(X): (P-\pi(z))f=0, f\sim 
\sum_{j\not\in  \tilde{\mcj}_{z_0}} c_j e^{i\tilde{r}_j(z)t}\phi_j\}.$$

To finish the proof, we need only understand what happens if 
$\Phi_j(w_J(z))$ has a pole at $z_0$ for some $j \in \tilde{\mcj}_{z_0}$.
However, a pole of $\Phi_j(w_J(z))$ cannot contribute to the singularity
of $R(z)$ at $z_0$ because the expansion on the ends of the singular part
is of the wrong form.  Therefore, the proof of this case follows much as
the proof of the previous one.
\end{proof}

The following proposition will also be useful.
\begin{prop}\label{p:notinspec}
Suppose $z_0\in \hat{Z}$ is such that $\pi(z_0)$ is not in the spectrum of $P$.
Then
$$m_{z_0}(R)=\mu d_{z_0}(S_{\mce_{z_0}}   ).$$
\end{prop}
\begin{proof}
We sketch the proof of this proposition, as it is very similar to the 
proof of Theorem \ref{thm:polesmult}.

Just as in Lemma \ref{l:mumsmaller}, one can show that 
$\mu m_{z_0}(S_{\mce_{z_0}}) \leq m_{z_0}(R)$.  Here, of
course, $$\dim \{f\in L^2(X): (P-\pi(z))f=0, f\sim 
\sum_{j\not\in  \tilde{\mcj}_{z_0}} c_j e^{ir_j(z)t}\phi_j\} =0,$$
 by our assumption on $\pi(z_0)$.  Moreover,
$S_{\mce_{z_0}}(z)$ has no zeros at $z_0$ as a zero would imply the existence 
of an $L^2$ eigenfunction, and thus
$$ \mu m_{z_0}(S_{\mce_{z_0}})= \mu d_{z_0}(S_{\mce_{z_0}}).$$

Finishing the proof thus requires showing that
$$ m_{z_0}(R)\leq \mu m_{z_0}(S_{\mce_{z_0}}).$$
This can be done as in Lemma \ref{l:mugreater}, first noting that, if
$\psi_l\not \equiv 0$ is in the image of the singular part of $R$ at $z_0$,
 then on the 
ends
\begin{equation}
\psi_l(t,y) =\sum_j \sum _{m\leq m_{j}} b_{ljm} t^m e^{i\tilde{r}_j(z_0)t}
\phi_j(y)
\end{equation} where $b_{ljm}\not =0$ for some $j\in \mce_{z_0}$ and some $m$.
In fact, because $\pi(z_0)$ is not an eigenvalue of $P$, the linear 
independence of $\psi_1$, $\psi_2$, ...$\psi_q$ in the image of the singular
part of $R(z)$ at $z_0$ is equivalent to the linear independence of
$$\sum_{j \in \tilde{\mce}_{z_0}} \sum _{m\leq m_{j}} 
b_{ljm} t^m e^{i\tilde{r}_j(z_0)t}
\phi_j(y),$$
$l=1,...,q$.
Using this and 
the fact that $\Xi_{z_0}(R)$ is symmetric, where $\Xi_{z_0}$ is the singular
part at $z_0$, we obtain that
\begin{multline*}
\dim \Image \Xi_{z_0}(R)  \\ \leq 
\dim \Image \Xi_{z_0}\left( 
\sum_{j,m\in \tilde{\mce}_{z_0}}S_{mj}(z) \Phi_m(w_J(z))\otimes 
\Phi_j(w_J(z))\frac{1}{\tilde{r}_j(z)}\right).
\end{multline*}
Then, if $\Phi_j(w_J(z))$ is regular at $z_0$ for all $j\in \tilde{J}_{z_0}$,
it is easy to see that 
$$m_{z_0}(R)\leq \mu m_{z_0} (S_{\mce_{z_0}}).$$
Again, if $\Phi_j(w_J(z))$ has a pole at $z_0$, it does not contribute to the 
singularities of $R(z)$ at $z_0$ but corresponds instead to a singularity
of $R(z)$ at $w_J(z_0)$.
\end{proof}

\section{Proof of Theorem \ref{thm:bbn}}

In this section we bound the number of poles of the resolvent in
neighborhoods of the 
ramification points on the boundary of the physical sheet, with
the size of the neighborhoods increasing.  We will use the fact that
if $\chi=1$ for $t<\max(a_i,0)$, then
$\dim \Image \Xi_{z_0}(R)= \dim \Image \Xi_{z_0}(R\chi)$, where
$\Xi_{z_0}(T)$ is the singular part of $T$ at $z_0$.

We recall Theorem \ref{thm:ln}
\begin{thmempty} 
Assume $X$ satisfies the hypothesis $(H1)$ and let $\beta <1$.
Then, in the connected components of $\{z'\in \hat{Z}: |r_m(z')|<\beta 
\sqrt{\alpha \nu_m}\}$ that
 meet the physical sheet, there are at most
${\mathcal O}(
m^{n-1})$ resonances.
\end{thmempty}
A corollary to this is

\begin{thm}\label{thm:bbn}
Assume $X$ satisfies the hypothesis $(H1)$ and let $\beta < 1$.
Then, on the closure of
the sheet with $\Im r_j(z)<0$ if and only if $j \leq m$, there are at
most ${\mathcal O}_{\beta}(m^{n-1})$ poles of the resolvent of $P$
 with $|r_m(z)|<\beta\sqrt{\alpha \nu_m}$.
\end{thm}

We shall use the Fredholm determinant method used in,
for example, \cite{gu-zwncrs,
pbdp, vodev, 
zwrp, zwspb}. 
We will find a trace class operator $K_m^{\#}(z)$ so that, in
the desired region, the 
poles of the resolvent are contained in the zeros of $I+ K_m^{\#}(z)$,
and thus in the zeros of $\det(I+ K_m^{\#}(z))$.  In addition, $K_m^{\#}(z_0)
=0$, where $z_0\in \hat{Z}$ is a ``base point'' which depends on $m$.  
To do this, we first construct 
an approximation of the resolvent adapted to this problem and valid in this
region, obtaining a compact operator
$K_m(z)$ so that the poles of the resolvent
in this region are 
contained in the zeros of $I+K_m(z)$.  Further manipulations simplify
$I+K_m(z)$.

For $m>1$, 
there are two ramification points on the boundary of the physical sheet
of $\hat{Z}$ which correspond to $\nu_m^2$,
and thus, for large $m$, two connected components 
of
$\{z'\in \hat{Z}: |r_m(z')|<\beta 
\sqrt{\alpha \nu_m}\}$. 
 We shall work near the ramification point
which is reached by taking a limit as $\Im \pi(z) \downarrow 0$, $\Re 
\pi(z )
\rightarrow \nu_m^2$.  We will designate this
ramification point by $(\nu_m^2)_+$.  
 A similar analysis works near the other point, $(\nu_m^2)_-$,
obtained by taking a limit for $z$ in the physical space, $\Im 
\pi(z) \uparrow 0$,
$\Re \pi(z) \rightarrow \nu_m^2$.

We assume $m>1$.

For $c>\max(a_i, 0)$, let 
$$X_c=X_{\comp}\sqcup _{i=1}^{m}(X_i\cap (a_i,c]
\times Y_i).$$
For $\zeta \in \Complex$, let $R_c(\zeta) =(P_{\rest X_c}-\zeta)^{-1}$
be the resolvent of $P$ on $X_c$ with Dirichlet boundary conditions.
(If $X$ has a boundary and we are considering Neumann boundary
conditions on $X$, we use Neumann boundary conditions on $X_c$ here.)
Recall that $R_Y(z)=(D_t^2+\Delta_Y-z)^{-1}$ is defined by (\ref{eq:RY}).

For $i=1,2,3$, choose $\chi_i\in C_c^{\infty}(X)$ so that 
$\chi_i\chi_{i+1}=\chi_i$, $i=1,2$, $\chi_1\equiv 1$ on 
$X_{\max(a_i,0)}$, $\nabla \chi_i$ only depends on $t$, and 
$|\nabla \chi_i |<\gamma$, $|\Delta \chi_i|< \gamma$.  It suffices to take
$\gamma \leq \alpha(1-\beta)^2\left(1000(1+\alpha)\right)^{-1}$.  
Choose $c_0$ so that the support of $\chi_3$
is properly contained in $X_{c_0}$.

Choose $z_0$ in the physical
plane so that $\pi(z_0)=\nu_m^2+\frac{1}{4}\alpha(1-\beta)^2\nu_mi$ and
let 
\begin{multline*}
E_m(z)=\chi_3 R_{c_0}(\pi(z))\Pi_m\chi_2 
+ \chi_3 R_{c_0}(\pi(z_0))(1-\Pi_m)\chi_2
+ (1-\chi_1)R_Y(z)(1-\chi_2).
\end{multline*}
Here $\Pi_M$ projects off of the eigenfunctions of 
$P_{\rest X_{c_0}}$ with eigenvalues in $(\nu_m^2-5\alpha \nu_m, 
\nu_m^2+5\alpha \nu_m)$.  This of course depends on $c_0$, but we omit this 
dependence in our notation.
Then
\begin{equation}(P-\pi(z))E_m(z)= I +\tilde{K}_m(z),
\end{equation}
where 
\begin{multline*}
\tilde{K}_m(z)= [P,\chi_3]
\left( R_{c_0}(\pi(z))\Pi_m+R_{c_0}(\pi(z_0))(1-\Pi_m)\right) \chi_2 \\
- \chi_3(\pi(z)-\pi(z_0))R_{c_0}(\pi(z_0))(1-\Pi_m)\chi_2 
-[P,\chi_1]R_Y(z)(1-\chi_2),
\end{multline*}
and, if the domain is restricted to 
$L^2_{\comp}(X)$, $\tilde{K}_m(z)$ is meromorphic on $\hat{Z}$,
with the poles corresponding to poles of 
$R_{c_0}(\pi(z))\Pi_m$.  Let $\chi \in C^{\infty}_c(X)$
be $1$ on the support of $\chi_3$.  Then
\begin{equation}
(P-\pi(z))E_m(z)\chi= \chi(I+\tilde{K}_m(z)\chi)=\chi(I+K_m(z))
\end{equation}
and $K_m(z)$ is compact.

Our choice of $z_0$ and $\gamma$ guarantee that $I+K_m(z_0)$ is 
invertible, with norm bounded by $2$.
  Then, by analytic Fredholm theory, $I+K_m(z)$ is 
invertible at all but a discrete set of points of $\hat{Z}$.  The 
points where it is not invertible correspond to zeros of 
$I+K_m(z)$, and the poles of the cut-off resolvent 
$(P-z)^{-1}\chi$ are included in the union of the zeros 
of $I+K_m(z)$ and the poles of $E_m(z)$.  However, we will restrict our
attention to a region that does not include any poles of $E_m(z)$.

Using the fact that $I+K_m(z_0)$ is invertible, with norm bounded by
$2$, we obtain that in the region in question,
the poles of the resolvent are contained in the zeros of 
\begin{equation}
I+(I+K_m(z_0))^{-1}(K_m(z)-K_m(z_0)).
\end{equation}

Now we restrict our attention to a region where
$|r_m(z)-r_m(z_0)|<\rho
\sqrt{\alpha\nu_m}$,
where 
$$ \rho = (\beta^2/4+3/4)^{1/2}$$
 and $z$ lies on one of the four
sheets that meet the ramification point $(\nu_m^2)_+$,
with $m$ large.  In this region,
we have
$$\|(I+K_m(z_0))^{-1}[P,\chi_3] 
\left(R_{c_0}(\pi(z))- R_{c_0}(\pi(z_0)\right)\Pi_m \chi_2 \|
\leq \frac{1}{3},$$
and 
\begin{equation}\label{eq:tdepenonly}
\|((I+K_m(z_0))^{-1}[P,\chi_1]
\sum_{\sigma_j^2>\nu_{m+1}^2} \left(R_{ej}(z)- R_{ej}(z_0)\right)(\chi-\chi_2)
\| \leq \frac{1}{3}.
\end{equation}
For (\ref{eq:tdepenonly}), we are using the fact that $[P,\chi_1]$ depends
only on $t$.
Therefore, in this region the poles of the resolvent are contained
in the zeros of 
\begin{equation}
I+K_{m}^{\#}(z),
\end{equation}
 where
\begin{equation*}
K_{m}^{\#}(z)=L_m(z)\left( K_{1m}(z)+K_{2m}(z)\right) 
\end{equation*}
with 
\begin{equation}
K_{1m}(z)=
-\chi_3(\pi(z)-\pi(z_0))R_{c_0}(\pi(z_0))(1-\Pi_m)\chi_2 ,
\end{equation}
\begin{equation}
K_{2m}(z)=-
[P,\chi_1]
\sum_{\sigma_j^2\leq \nu_{m+1}^2} 
\left( R_{ej}(z)- R_{ej}(z_0)\right)(\chi-\chi_2),
\end{equation}
and the norm of 
$L_m(z)$
 is bounded, independent of $m$ and $z$, as long
as we stay in the region described.

Now $K_m^{\#}(z)$ is trace class.  We consider the function
$$h(z)=\det(I+K^{\#}_m(z)),$$ and note that $h(z_0)=1$ and that $h$ is 
holomorphic 
on $\hat{Z}$
when $|r_m(z)-r_m(z_0)|<2\sqrt{\alpha \nu_m}$.  We will apply Jensen's 
theorem to $h$ to obtain an upper bound on the number of zeros of 
the resolvent in this region, and to do this we need the following 
lemma.

\begin{lemma}
If $|r_m(z)-r_m(z_0)|
\leq \rho \sqrt{\alpha \nu_m}$
 and $z$ lies in the connected component
of $\{z'\in \hat{Z}: |r_m(z')-r_m(z_0)|
\leq\rho \sqrt{\alpha \nu_m}\}$
 containing $(\nu_m^2)_+$, 
then
$|h(z)|\leq Ce^{C\langle m \rangle ^{n-1}}$.
\end{lemma}
\begin{proof}

We remark that because of Weyl's law and the hypothesis (H1), for large
$m$ $\nu_m$ is bounded above and below by constant multiples of $m$.

We use the property that 
\begin{equation}\label{eq:adddet}
|\det(I+A+B)|\leq \det(I+|A|)^2 \det(I+|B|)^2
\end{equation}
(e.g. \cite[Lemma 6.1]{gu-zwncrs}).  
Since $(I-\Pi_m)\chi_2$ has rank bounded by $C m^{n-1}$ and 
$\| R_{c_0}(\pi(z_0)) \| \leq \frac{C}{\langle \nu_m \rangle}$, we obtain,
in the region with $|r_j(z)-r_j(z_0)|\leq 
\rho
\sqrt{\alpha \nu_m}$, and thus $|\pi(z)-\pi(z_0)|\leq C\langle \nu_m\rangle$,
$$|\det(1+|L_m(z)K_{1m}(z)|)| \leq Ce^{C\langle m \rangle^{n-1} }.$$

Consider next 
\begin{multline*}
L_m(z)K_{2m}(z)=L_{m}(z) [P,\chi_1]
\sum_{\sigma_j^2\leq \nu_{m+1}^2} 
\left( R_{ej}(z)- R_{ej}(z_0)\right)(\chi-\chi_2)\\
 = L_m(z)[P,\chi_1]
\sum_{j=1}^m
\left( e^{i|t-t'|r_j(z)}
-e^{i|t+t'|r_j(z)}
-e^{i|t-t'|r_j(z_0)} +e^{i|t+t'|r_j(z_0)}\right) \\ 
\times \sum_{\sigma_l^2=\nu_j^2}
\phi_l(y)\overline{\phi_l}(y')(\chi-\chi_2).
\end{multline*}
Since $[P,\chi_1]$ and $\chi-\chi_2$ have disjoint supports, 
$$
[P,\chi_1]
\left( e^{i|t-t'|r_j(z)}
-e^{i|t+t'|r_j(z)}
-e^{i|t-t'|r_j(z_0)} +e^{i|t+t'|r_j(z_0)}\right)
\phi_l(y)\overline{\phi_l}(y')(\chi-\chi_2)
$$ 
is a rank four operator with norm bounded by 
$Ce^{C|\Im r_j(z)|}$.
We need, therefore, to bound
\begin{equation}\label{eq:imsum}
\sum_1^{m+1} |\Im r_j(z)|\# \{\sigma_l^2:\sigma_l^2=\nu_j^2\}.
\end{equation}
We have
$$\# \{\sigma_l^2:\sigma_l^2=\nu_j^2\} \leq C j^{n-2}.$$
On the region of interest, $|\Im r_j(z)|$ can be bounded by 
$2\alpha \nu_m(|\nu_m^2-\nu_j^2|)^{-1/2}$, if $j< m$, and $\sqrt{\alpha 
\nu_m}$, if $j= m$.
Therefore,
\begin{equation}
\sum_1^{m+1} |\Im r_j(z)|\# \{\sigma_l^2:\sigma_l^2=\nu_j^2\}
= \sum_1^{m-2} |\Im r_j(z)|\# \{\sigma_l^2:\sigma_l^2=\nu_j^2\}+ {\mathcal O}
(m^{n-1}).
\end{equation}
To bound (\ref{eq:imsum}), then, it suffices to bound
\begin{multline}\label{eq:imsum2}
 2\alpha \nu_m \int_0^{\nu_{m-1}}
\left( \nu_m^2-l^2\right)^{-1/2}
l^{n-2} dl   \\
= 2\alpha \nu_m^{n-1}\int _0^{\nu_{m-1}/\nu_m}
(1-s^2)^{-1/2}s^{n-2} ds
\leq \alpha C \nu_m^{n-1}.
\end{multline}
\end{proof}

\begin{proof}
[Proof of Theorem \ref{thm:ln}]
The poles of the resolvent in the region in question are 
contained in the zeros of $h(z)$.  Working near the point $(\nu_m^2)_+$,
we use $r_m(z)$ as a coordinate, which we may do as long as we keep 
away from other ramification points or regions where $r_m(z)$ fails to
be one-to-one.
We apply Jensen's Theorem to $h$, 
using a circle centered at $z_0$ and having radius 
$\rho
\sqrt{\alpha \nu_m}$.
Then the theorem follows because the disk $|r_m(z)|\leq \beta \sqrt{\alpha 
\nu_m}$ is properly contained in the disk $|r_m(z)-r_m(z_0)| 
\leq \rho
 \sqrt{\alpha \nu_m}$, with the ratio of the distance between the two 
boundaries
 and the radius of the larger disk bounded from below, independent
of $m$.
\end{proof}

\section{Proof of Theorem \ref{thm:ppp}}

We recall the statement of Theorem \ref{thm:ppp}:
\begin{thmempty}
Let $X= \Real \times Y$ and suppose $X$ satisfies (H1), and let 
$\rho>0$ be fixed.  Consider the operator $\Delta +V$,
for real-valued $V\in L^{\infty}_{\comp}(X)$.  
Then, on the connected components of $\{z\in \hat{Z}: |r_m(z)|<\rho\}$
that meet the physical sheet, the number of poles is bounded by
$C(1+m^{n-2})$.
\end{thmempty}

We prove this by the Fredholm determinant method as in the previous theorem.
We assume that $\rho>1$.

Let $R_0(z)=(\Delta -z)^{-1}$; its Schwartz kernel is given by
\begin{equation}
\label{eq:0resolvent}
R_0(z)=\sum_{j=1}^{\infty}\frac{i}{2r_j(z)}e^{i|t-t'|r_j(z)}
\sum_{\sigma_l^2=\nu_j^2}\phi_l(y)\overline{\phi}_l(y'),
\end{equation}
and let $R_V(z)=(\Delta +V-z)^{-1}$.  We have 
\begin{equation}\label{eq:potr}
(\Delta +V-\pi(z))R_0(z)=I+VR_0(z).
\end{equation}
Let $\chi(t)\in C_c^{\infty}(X)$ be one on the support of $V$,
with $|\chi|\leq 1$.
Then, if $z_0\in \hat{Z}$ is not a 
ramification point  and $R_V(z)\chi$ has a pole at $z_0$, then
$I+VR_0(z_0)\chi$ has nontrivial null space.  

For $r>0$, $m\in \Natural$, set $B_{m,r}$ to be the connected components
of $\{ z \in \hat{Z}: |r_m(z)|<r\}$ that meet the physical sheet.

Let $\tilde{\rho}=\max(\rho, 2 \|V\|_{\infty}^{1/2})$, and restrict $z$ to 
$B_{m,4\tilde{\rho}}$.  Take $m$ sufficiently large that the only 
ramification points in $B_{m,4\tilde{\rho}}$ correspond to $\nu_m^2$.
Then 
$$I+VR_0(z)\chi= (I+K_1(z))(I+(I+K_1(z))^{-1}K_2(z))$$
where $K_1(z)$ has Schwartz kernel
$$K_1(z)= V(t,y)\sum_{j\not = m}\frac
{i}{2r_j(z)}e^{i|t-t'|r_j(z)}
\sum_{\sigma_l^2=\nu_j^2}\phi_l(y)\overline{\phi}_l(y')\chi(t')$$
and $K_2(z)$
$$K_2(z) = V(t,y)\frac
{i}{2r_m(z)}e^{i|t-t'|r_m(z)}
\sum_{\sigma_l^2=\nu_m^2}\phi_l(y)\overline{\phi}_l(y')\chi(t').$$
By choosing $m$ sufficiently large, so that 
$$\frac{A^{1/2}}{|\alpha \nu_m-16\tilde{\rho}^2|^{1/2}}
\exp\left( A \frac{16 \tilde{\rho}^2}{\sqrt{2\alpha \nu_m}}\right) 
< \frac{1}{2\|V\|},$$
where $A=\max_{t,t'\in \supp \chi}|t-t'|$,
we have 
$\|(I+K_1(z))^{-1}\|\leq 2$ on $B_{m, 4 \tilde{\rho}}$.
The poles of $R_V(z)$ in $B_{m,4\tilde{\rho}}$, other than the ramification 
point, correspond to values of 
$z$ for which $I+(I+K_1(z))^{-1}K_2(z)$ has non-trivial null space.  We
remark that $R_0(z)$ has a pole of rank $M_Y(\nu_m^2)$ at the ramification
point, where $M_Y(\nu_m^2)$ is the multiplicity of $\nu_m^2$ as an 
eigenvalue of $\Delta_Y$, and this can contribute a pole of up to the same
multiplicity
 to $R_V(z)$, even if $I+(I+K_1(z))^{-1}K_2(z)$ is invertible there.

Let 
$$K_3(z)=(I+K_1(z))^{-1}K_2(z).$$  
The poles of $R_V(z)$ in $B_{m,4\tilde{\rho}}$ are contained in the 
zeros of $I+K_3(z)$ in the same region, except, possibly,
at the ramification points, as discussed above.  Now we shall work
 near $(\nu_m^2)_+$ as in the
proof of the previous theorem, as a similar analysis will work for
the other connected component of $B_{m, 4\tilde{\rho}}$.
Choose $z_0$ in the physical space
with
$\pi(z_0)=\nu_m^2+4i\|V\|_{\infty}$.  Then
$\|K_2(z_0)\|\leq 1/4$, and 
$$\left\|\left(I+(I+K_1(z_0))^{-1}K_2(z_0)\right)^{-1}\right\|
= \|(I+K_3(z_0))^{-1}\|\leq 2.$$
Let 
$$h(z)=\left( \frac{r_m(z)}{r_m(z_0)}\right)^{2M_Y(\nu_m^2)}
\det\left( I + (I+K_3(z_0))^{-1}(K_3(z)-K_3(z_0))\right).$$
Then $h(z_0)=1$, and, except, possibly, for some at the 
ramification point, the poles of $R_V(z)$ in the connected 
component of $B_{m,4\tilde{\rho}}$ that includes $(\nu_m^2)_+$,
are contained in the zeros of $h(z)$ in the same region.  This 
misses at most $M_Y(\nu_m^2)={\mathcal O}(\nu^{n-2})$ poles of
$R_V$ at the ramification point.  An application of Jensen's theorem
on a circle centered at $z_0$ and with $|r_m(z)-r_m(z_0)|\leq 3\tilde{\rho}$
will then finish the proof, after we have proved the following lemma.
\begin{lemma}
On the connected component of $B_{m,4\tilde{\rho}}$ that contains
$(\nu_m^2)_+$, 
$$|h(z)|\leq \exp(C(1+|m|^{n-2})).$$
\end{lemma}

\begin{proof}
Let
$R_{ems}(z)$ be the operator with Schwartz kernel
$$\frac{i}{2r_m(z)}\chi(t)
(e^{i|t-t'|r_m(z)}-1)
$$
and let 
$R_{em0}(z_0)$ be the operator with Schwartz kernel
$$\frac{i}{2r_m(z_0)}\chi(t) e^{i|t-t'|r_m(z_0)}.$$
Let 
$$A_1(z)=(I+K_1(z))^{-1}V\frac{i}{2r_m(z)}
\sum_{\sigma_l^2=\nu_m^2}
\phi_l\otimes \overline{\phi}_l\chi.$$
Then
\begin{multline*}
K_3(z)-K_3(z_0) =\left(
(I+K_1(z))^{-1}VR_{ems}(z)-(I+K_1(z_0))^{-1}VR_{em0}(z_0)\right)\\ \times
\sum _{\sigma_l^2=\nu_m^2}
\phi_l\otimes \overline{\phi}_l\chi +A_1.
\end{multline*}
Now we shall use (\ref{eq:adddet}) and 
$$|\det(I+|BT|)| \leq \det(I+\|B\| |T|)
$$
(e.g. \cite[Lemma 6.1]{gu-zwncrs}).
Then
\begin{multline*}
|\det(I+(I+K_3(z_0))^{-1}(K_3(z)-K_3(z_0)))| \\
\leq \det 
(I+\|(I+K_3(z_0))^{-1}(I+K_1(z))^{-1}V\| |R_{ems}(z)\sum_{\sigma_l^2=\nu_m^2}
\phi_l\otimes \overline{\phi_l} \chi|)^4
\\ \times 
\det(I+\|(I+K_3(z_0))^{-1}(I+K_1(z_0))^{-1}V\| |R_{em0}(z_0)\sum_{\sigma_l^2=\nu_m^2}\phi_l\otimes \overline{\phi_l} \chi|)^4 \\ \times 
\det(I+|(I+K_3(z_0))^{-1}A_1|)^2
.
\end{multline*}

On $B_{m,4\tilde{\rho}}$, $\|(I+K_3(z_0))^{-1}(I+K_1(z))^{-1}\| \leq C$.
  For a 
compact operator $A$, let $\mu_1(A)\geq \mu_2(A)\geq \mu_3(A)\geq ...$
be the 
characteristic values of $A$; that is, the eigenvalues
of $|A^*A|^{1/2}$.  Then, if $A$ is trace class, 
$$\det(I+|A|)=\prod   (I+\mu_j(A)).$$
For $p=0, 1,2,...$, 
$$\mu_{pM_Y(\nu_m^2)+j}(R_{ems}(z)\sum_{\sigma_l^2=\nu_m^2} \phi_l
\otimes
\overline{\phi_l}\chi)=\tilde{\mu}_p(r_m(z)),\; j=1,...,M_Y(\nu_m^2).$$
Here we use the fact that $R_{ems}(z)$ depends on $z$ only through $r_m(z)$.
The $\tilde{\mu}_p(w)$ are independent of $m$.
Therefore,
\begin{multline*}
|\det (I+\|(I+K_3(z_0))^{-1}(I+K_1(z))^{-1}V\| | R_{ems}(z)
\sum_{\sigma_l^2=\nu_m^2}
\phi_l\otimes
\overline{\phi_l} \chi|)| \\ \leq C^{M_Y(\nu_m^2)}\leq C^{\langle m \rangle
^{n-2}}.
\end{multline*} The same argument gives
$$|\det (I+\|(I+K_3(z_0))^{-1}(I+K_1(z_0))^{-1}V\| | R_{em0}(z_0)
\sum_{\sigma_l^2=\nu_m^2}
\phi_l\otimes
\overline{\phi_l} \chi|)| \leq C^{\langle m \rangle
^{n-2}}.$$

A similar argument bounds 
$$r_m(z)^{2M_Y(\nu_m^2)}
\det(I+|(I+K_3(z_0))^{-1}A_1|)^2.$$  The rank of 
$(I+K_3(z_0))^{-1}A_1$ is $M_Y(\nu_m^2)$, and 
$\|(I+K_3(z_0))^{-1}A_1\|$
is bounded for
$4\tilde{\rho}\geq |r_m(z)|>c>0$.  Since 
$r_m(z) (I+K_3(z_0))^{-1}A_1$ is bounded on $B_{m,4\tilde{\rho}}$,
$$r_m^{M_Y(\nu_m^2)}\det(I+|(I+K_3(z_0))^{-1}A_1|)\leq C^{M_Y(\nu_m^2)}
\leq C^{\langle m \rangle ^{n-2}}$$
on $B_{m,4\tilde{\rho}}$.
\end{proof}

A consequence of Theorem \ref{thm:ppp} is
\begin{cor}
Let $X=\Real \times Y$ and suppose $X$ satisfies (H1).  Let
$V\in L^{\infty}_{\comp}(X)$ be real-valued, and let 
$$N(\lambda)=\# \{\lambda_j^2\leq \lambda^2: \lambda_j^2 \;
\text{ is an eigenvalue of } \Delta +V\}.$$
Then
$$N(\lambda) ={\mathcal O}(\lambda^{n-1}).$$
\end{cor}
\begin{proof}
Suppose $\tau \in \Real_+$.  Then, using (\ref{eq:0resolvent}) and
(\ref{eq:potr}), we se that if
$$\frac{\|V\| (\length \supp V + 1)}{\min_{j}|r_j(\tau)|}\leq \frac
{1}{2}$$
then $\tau $ cannot be an eigenvalue of $\Delta +V$.  Here 
$$\length \supp V= \max_{a,b\in \supp V}|a-b|.$$  This means that any
eigenvalue must lie within a fixed distance of some $\nu_j^2$.  Theorem
\ref{thm:ppp} provides a bound on the number of eigenvalues within such
a ball; summing over the $\nu_j^2$ we obtain the corollary.
\end{proof}

\section{Proof of Theorem \ref{thm:fixedsheet}}

We recall Theorem \ref{thm:fixedsheet}:
\begin{thmempty}
Fix a sheet of $\hat{Z}$, and let $\{z_k\}$ be the resonances
of $P$ on this 
sheet.  Then 
$$\sum 
\frac{|\Im (r_1(z_k))|}{|r_1(z_k)|^2}< \infty.$$
\end{thmempty}

In the proof of this theorem we shall use Proposition \ref{p:notinspec}
and Carleman's Theorem, which
we recall (e.g. \cite[Section V.i]{levin}).
\begin{thmempty}[Carleman]
If $F(\zeta)$ is a holomorphic function in the region $\Im \zeta \geq 0$,
$F(0)=1$, 
and if $a_k=r_ke^{i\theta_k}$ ($k=1,2,...$) are its zeros in this 
region, then
\begin{multline*}
\sum _{r_k \leq R}(\frac{1}{r_k}-\frac{r_k}{R^2}) \sin \theta_k 
= \frac{1}{\pi R} \int _0^{\pi}
\ln |F(Re^{i\theta})|\sin \theta d \theta  \\
+ \frac{1}{2\pi} \int_0^R \left(\frac{1}{x^2} -\frac{1}{R^2}\right)
\ln |F(x)F(-x)|dx + \frac{1}{2}\Im F'(0).
\end{multline*}
\end{thmempty}
We note that Carleman's Theorem also holds for a function $F(\zeta)$ which
is holomorphic in $\Im \zeta >0$ and continuous in $\Im \zeta \geq 0$.  
In order to see this, apply Carleman's Theorem to $F_{\epsilon}(\zeta)
= F(\zeta+i\epsilon)/F(i\epsilon)$, $\epsilon >0$.  
Then, since both sides of the 
equation are continuous in $\epsilon$ for small $\epsilon \geq 0$, the theorem
holds in this case as well.

\begin{proof}[Proof of Theorem \ref{thm:fixedsheet}]

Fix a sheet of $\hat{Z}$, and let 
$$\mce= \{ j\in \Natural: \Im r_j(z)<0 \; \text{ on this sheet}\}.$$
By Proposition \ref{p:notinspec},
 the poles of the resolvent on this sheet (but not on
its boundary) correspond, with multiplicity, to the poles of 
$\det S_{\mce}(z)$ on this sheet.
We have, by Proposition \ref{prop:sinverse},
$$S_{\mce}^{-1}(z)=S_{\mce}(w_{\mce}(z)),$$
and, if $z$ lies on the sheet with
$$\Im r_j(z)<0 \text{ if and only if}\; j\in \mce,$$
then
 $w_{\mce}(z)$ lies on the physical sheet.
Therefore, we reduce the problem to a question about the zeros of 
$\det S_{\mce}(z)$ for $z$ on the physical sheet.

It is helpful to identify the physical sheet with the upper half plane 
using the variable $\zeta =r_1(z)$.  Let $\Psi_{ij}(\zeta)=S_{ij}(z(\zeta))$
and $\Psi(\zeta)=
S_{\mce}(z(\zeta))$.  Using the fact that $S_{\mce}(z)$ is meromorphic on 
$\hat{Z}$ we can extend $\Psi$ to the closed upper half plane by
continuity, except,
perhaps, for a finite number of points corresponding to poles of 
$S_{\mce}(z)$.  We shall also call these points poles of $\Psi(\zeta)$.
The matrix $\Psi(\zeta)$ has at most a finite number of poles in the 
closed upper half plane.

To prove the theorem, we shall use Carleman's theorem applied to a multiple
of $\det \Psi(\zeta)$, chosen so that the product is 
holomorphic in the upper half plane and continuous on its closure. 
 In order to do this, we need bounds on $\det 
\Psi$.  Let $J=\max\{ j\in \mce\}$.
If $\zeta \in \Real$, $\zeta^2+\sigma_1^2>\nu_J^2$, then 
$|\Psi_{ij}|\leq C$ (\cite{ace}).  To bound $|\Psi_{ij}(\zeta)|$
away from the real axis, we need to bound $|S_{ij}(z)|$, where 
$S_{ij}$ is determined by the expansion on the ends of $\Phi_j$,
see (\ref{eq:phij1}) and (\ref{eq:phij2}).
We recall that
\begin{equation}\label{eq:phiconst}
\Phi_j(z(\zeta))=\chi(t)e^{-i\tilde{r}_j(\zeta^2+\sigma_1^2)t}\phi_j
-(P-\zeta^2-\sigma_1^2)^{-1}(P-\zeta^2-\sigma_1^2)
\chi(t)e^{-i\tilde{r}_j(\zeta^2+\sigma_1^2)t}\phi_j
\end{equation}
where $\chi(t)\in C^{\infty}(\Real)$ is supported in $t>\max(a_i)$
and is one in a neighborhood of infinity.

We note that for $j\in \tilde{\mce}$, 
\begin{multline}\label{eq:normest}
\| (P-\zeta^2 -\sigma_1^2) 
\chi(t)e^{-i\tilde{r}_j(\zeta^2+\sigma_1^2)t}\phi_j \|
\leq C \langle \zeta^2+\sigma_1^2-\sigma_j^2\rangle ^{1/2}e^{C |\Im 
\tilde{r}_j(\zeta^2+\sigma_1^2)|}\\
\leq  C \langle \zeta^2+\sigma_1^2-\sigma_j^2\rangle ^{1/2}e^{C \langle 
\Im \zeta\rangle }
\end{multline}
and
\begin{equation}\label{eq:minest}
\| \left(e^{i\tilde{r}_l( \zeta^2+\sigma_1^2)t}\phi_l\right)
_{\rest t>a}\|
= (2|\Im\tilde{r}_l(\zeta^2+\sigma_1^2)|)^{-1/2}
 e^{-a\Im \tilde{r}_l(\zeta^2+\sigma_1^2)}.
\end{equation}

Since 
$$\| (P-\zeta^2-\sigma_1^2)^{-1}\| 
\leq (\dist(\zeta^2+\sigma_1^2, \sigma(P))^{-1},$$
we have, using (\ref{eq:phij2}), (\ref{eq:phiconst}), (\ref{eq:normest}) and (\ref{eq:minest}),
$$|\Psi_{lj}(\zeta)|\leq 
\frac{C\langle \zeta^2+\sigma_1^2-\sigma_l^2\rangle ^{1/2}}
{\dist(\zeta^2+\sigma_1^2, \sigma(P))}e^{C \langle \Im \zeta \rangle}$$
when $|\zeta|$ is large.  This proves
$$|\det \Psi(\zeta)|\leq Ce^{C\langle \Im \zeta \rangle}$$
if
$\Im \zeta>\epsilon >0$ and $|\zeta|$ large.  To obtain a bound in the 
closure of the upper half-plane, apply the Phragmen-Lindel\"of theorem
to 
$$h(\zeta)=\prod_{1}^{k_0}\frac{\zeta -\zeta_j}{\zeta -\zeta_j+Mi}\det 
\Psi(\zeta)$$
where $\zeta_1,$ $\zeta_2$, ... $\zeta_{k_0}$ are the poles of $\det
\Psi(\zeta)$ 
in the closed upper half plane and $M$ is chosen sufficiently large
that $M>\Im k_j$, $j=1,$...,$k_0$.  Then $h(\zeta)$ is holomorphic in the
upper half plane, continuous in the closed upper half plane, and 
$$|h(\zeta)|\leq Ce^{C\langle \Im \zeta\rangle}$$
in the closed upper half plane.  
An application of Carleman's Theorem to $h(\zeta)$ finishes the proof.
\end{proof}

\small
\noindent
{\sc 
Department of Mathematics,
University of Missouri
Columbia, Missouri 65211\\}
\end{document}